\documentclass{article}
\usepackage{graphicx} 
\usepackage{amsbsy}
\usepackage{bm}
\usepackage{amssymb,amsmath,amsfonts,amsthm}
\usepackage{color}
\newtheorem {Proposition}{Proposition}[section]
\newtheorem {Lemma}{Lemma}[section]

\newtheorem {corollary}{Corollary}[section]

\newtheorem {Remark}{Remark}[section]
\usepackage{amssymb,amsmath}
\usepackage{authblk}
\usepackage{epsfig}
\usepackage{enumerate,graphicx}
\usepackage{amsfonts}
\usepackage{graphicx}
\DeclareGraphicsExtensions{.pdf, .jpg, .png , .bmp, .ps}
\def\build#1_#2^#3{\mathrel{\mathop{\kern 0pt#1}\limits_{#2}^{#3}}}

\title{Some remarks on Gordin-Lif\v{s}ic's condition for 
martingale approximations}
\author[1]{J\`er\^ome Dedecker}
\affil[1]{Universit\'e  Paris Cit\'e,  Laboratoire MAP5 and CNRS UMR 8145}
\author[2]{Florence Merlev\`ede}
\affil[2]{Université Gustave Eiffel, LAMA and CNRS UMR 8050}
\date{}

\begin{document}

\maketitle

\begin{abstract}
In this note, we study a condition introduced by Gordin and Lif\v{s}ic in 1981 to establish the Central Limit Theorem for additive functionals of stationary Markov chains with normal transition operator.  In the more general setting of strictly stationary sequences satisfying the Gordin-Lif\v{s}ic condition, we give sufficient (and sometimes also necessary) conditions for partial sums to be approximated in ${\mathbb L}^2$ by a martingale with stationary increments. Various types of ${\mathbb L}^2$ approximations are described, leading to different versions of the central limit theorem (annealed, quenched, functional form...). The optimality of the conditions is discussed, and an application to the class of semi-linear processes is presented.
\end{abstract}


\section{Introduction}
Let $(\Omega, {\mathcal A}, {\mathbb P})$ be a probability space and let $T : \Omega \rightarrow \Omega $ be a bijective bimeasurable transformation that preserves the probability ${\mathbb P}$, and let ${\mathcal I}$ be the sub-$\sigma$-algebra of ${\mathcal A}$ consisting of all invariant sets by $T$. Let ${\mathcal F}_0$ be a $\sigma$-algebra of ${\mathcal A}$ satisfying ${\mathcal F}_0 \subseteq T^{-1}({\mathcal F}_0)$, and let ${\mathcal F}_i = T^{-i}({\mathcal F}_0)$. Let $X_0$ be a ${\mathcal F}_0$-measurable, square integrable and centered random variable, and define the stationary sequence $(X_i)_{i \in {\mathbb Z}}$ by $X_i=X_0 \circ T^i$. Let then 
$$
S_n = X_1 + X_2 + \cdots + X_n\, , \quad \text{and for $t \in [0,1]$,} \quad W_n(t)=S_{[nt]} + (nt-[nt])X_{[nt]+1}\, .
$$

The sequence $(X_i, {\mathcal F}_i)_{i \geq 0}$ is said to satisfy the Gordin-Lif\v{s}ic condition if
\[
\sum_{k \geq 0} \|{\mathbb E}(X_k | {\mathcal F}_0)\|_2^2 < \infty \, .  \tag{G-L}  \label{GL}
\]
This condition was introduced by Gordin and Lif\v{s}ic \cite{GL81} in the case where $X_k=f(\xi_k)$ is a function of a stationary Markov chain $(\xi_k)_{k \in {\mathbb Z}}$ and ${\mathcal F}_k=\sigma( \xi_i, i \leq k)$. They showed that if the sequence $(\xi_i)_{i \geq 0}$ is ergodic and if the transition operator $K$ of the chain is normal (that is $K^*K=KK^*$, where $K^*$ is the adjoint operator of $K$), then $n^{-1/2}S_n$ satisfies the Central Limit Theorem (CLT). We refer to the book by Borodin and Ibragimov \cite{BI94} or to the paper by Derrienic and Lin \cite{DL01} for a proof of this remarkable result. Condition \eqref{GL} has also been investigated in detail by Cuny and Lin \cite{CL15}: see Proposition 7.1, Proposition 7.4, Theorem 8.1 and Theorem 8.2 there. Still in the context where $X_k=f(\xi_k)$, Peligrad \cite{P24} recently proved that $n^{-1/2}S_n$ satisfied the quenched CLT provided that $\sum_{k \geq 0} \|{\mathbb E}(X_0 | \xi_k, \xi_{-k})\|_2^2 < \infty$, which is a reinforcement of \eqref{GL}.

Dedecker \cite{D15} showed that there exists a stationary ergodic sequence $(X_i)_{i \geq 0}$ (in fact a linear process with independent and identically distributed (iid) innovations) satisfying 
\begin{equation}\label{GLlog}
\sum_{k \geq 1}  \log (k)\|{\mathbb E}(X_k | {\mathcal F}_0)\|_2^2 < \infty \, ,
\end{equation}
but $n^{-1/2}S_n$ is not stochastically bounded (see also Proposition \ref{CE} of Section \ref{QLin}). This proves that the condition \eqref{GL} alone is not sufficient for the CLT with normalization $\sqrt n$.

\medskip

In this paper, we consider the following question: if  \eqref{GL} holds, under what additional condition do we get a martingale approximation in ${\mathbb L}^2$? More precisely, we will consider different kinds of approximation, which we now describe. 

\medskip

\noindent $\bullet$ We say that there is a martingale approximation in ${\mathbb L}^2$ (denoted by ${\mathbb L}^2$-MA) if there exists a (necessarily unique) ${\mathcal F}_0$-measurable random variable $D$ in ${\mathbb L}^2$ such that ${\mathbb E}(D \circ T|{\mathcal F_0})= 0 $ a.s. and 
$$
\lim_{n \rightarrow \infty} \frac 1 n \left \| S_n - \sum_{k=1}^n D\circ T^k \right \|_2^2 = 0 \, .
$$
A well known consequence of  ${\mathbb L}^2$-MA is that $n^{-1/2} S_n $ converges in distribution to $\sqrt{{\mathbb E}(D^2|{\mathcal I})} Z$, where $Z$  is ${\mathcal N}(0,1)$-distributed and independent of ${\mathcal I}$.

\medskip

\noindent $\bullet$
We say that there is a maximal martingale approximation in ${\mathbb L}^2$ (denoted by ${\mathbb L}^{2}$-MMA) if there exists a (necessarily unique) ${\mathcal F}_0$-measurable random variable $D$ in ${\mathbb L}^2$ such that ${\mathbb E}(D \circ T|{\mathcal F_0})= 0 $ a.s. and 
$$
\lim_{n \rightarrow \infty} \frac 1 n \left \|\max_{1 \leq k \leq n} \left | S_k - \sum_{i=1}^k D\circ T^i \right | \right \|_2^2 = 0 \, .
$$
A well known consequence of  ${\mathbb L}^2$-MMA is that the Weak Invariance Principle (WIP) holds, that is :  $n^{-1/2} W_n $ converges in distribution in $(C([0,1], \|\cdot \|_\infty)$ to $\sqrt{{\mathbb E}(D^2|{\mathcal I})} W$, where $W$  is a standard Wiener process independent of ${\mathcal I}$.

\medskip

\noindent $\bullet$
We say that there is a martingale approximation in ${\mathbb L}^2({\mathcal F}_0)$ (denoted by ${\mathbb L}^2_0$-MA) if there exists a (necessarily unique) ${\mathcal F}_0$-measurable random variable $D$ in ${\mathbb L}^2$ such that ${\mathbb E}(D \circ T|{\mathcal F_0})= 0 $ a.s. and 
$$
\lim_{n \rightarrow \infty} \frac 1 n {\mathbb E} \left ( \left ( S_n - \sum_{k=1}^n D\circ T^k \right )^2 \Big | {\mathcal F}_0 \right ) = 0 \, , \, \text{almost surely and in ${\mathbb L}^1$.}
$$
Assume moreover that there exists a regular version ${\mathbb P}_{T|{\mathcal F}_0}$, and that the conditional expectations of $g\circ T$ with respect to ${\mathcal F}_0$  are taken through this conditional probability (we will always assume that this is the case when dealing with quenched results). An easy consequence of  ${\mathbb L}^2_0$-MA is that the quenched version of the CLT (Quenched CLT) holds, that is:  on a set of probability 1, for any continuous and bounded function $\varphi$, 
$$
\lim_{n \rightarrow \infty} {\mathbb E}(\varphi(n^{-1/2} S_n) |{\mathcal F_0}) = \int \varphi \left ( \sqrt{{\mathbb E}(D^2|{\mathcal I})}x \right )f(x) dx \, , $$
where $f$ is the density of a standard normal random variable. 

\medskip

\noindent $\bullet$
We say that there is a maximal martingale approximation in ${\mathbb L}^2({\mathcal F}_0)$ (denoted by ${\mathbb L}^2_0$-MMA) if there exists a (necessarily unique) ${\mathcal F}_0$-measurable random variable $D$ in ${\mathbb L}^2$ such that ${\mathbb E}(D \circ T|{\mathcal F_0})= 0 $ a.s. and 
$$
\lim_{n \rightarrow \infty} \frac 1 n {\mathbb E} \left (  \max_{1 \leq k \leq n} \left( S_k - \sum_{i=1}^k D\circ T^i \right )^2 \Big | {\mathcal F}_0 \right ) = 0 \, , \text{almost surely and in ${\mathbb L}^1$.}
$$
An easy consequence of  ${\mathbb L}^2_0$-MMA is that the quenched version of the WIP (Quenched WIP) holds, that is:  on a set of probability 1, for any continuous and bounded function $\varphi$, 
$$
\lim_{n \rightarrow \infty} {\mathbb E}(\varphi(n^{-1/2} W_n) |{\mathcal F_0}) = \int \varphi \left ( \sqrt{{\mathbb E}(D^2|{\mathcal I})}x \right )W(dx)  \, , $$
where $W$ is the distribution  of a standard Wiener process. 

\medskip

Note that we have the hierarchy: ${\mathbb L}^2_0$-MMA $\Rightarrow$ ${\mathbb L}^2_0$-MA $\Rightarrow$ ${\mathbb L}^2$-MA, and also ${\mathbb L}^2_0$-MMA $\Rightarrow$ ${\mathbb L}^2$-MMA.

\medskip As we already mentioned, our question is the following: if  \eqref{GL} holds, under what additional  condition do we get the martingale approximations described above? Let us first mention here that, in the case where $X_k=f(\xi_k)$ is a function of a stationary Markov chain with normal transition operator, then \eqref{GL} implies  ${\mathbb L}^2$-MA, as proved by Gordin and Peligrad \cite{GP11}, Remark 11.

Note that the strongest martingale approximation, ${\mathbb L}^2_0$-MMA,  holds provided 
\begin{equation}\label{MWstrong}
\sum_{k \geq 1} \frac{\|{\mathbb E}(X_k | {\mathcal F}_0)\|_2}{\sqrt k} < \infty \, .
\end{equation}
This has been proved by Cuny and Voln\'y \cite{CV13}, Corollary 2. It is also a direct consequence of Theorem 2.7 of Cuny and Merlev\`ede \cite{CM14}, who proved that ${\mathbb L}^2_0$-MMA holds under the Maxwell-Woodroofe condition \cite{MW00}. It is easy to see that \eqref{MWstrong} $\Rightarrow$ \eqref{GL} (see Proposition \ref{comp} of Section \ref{Sec4}). Both conditions \eqref{GL} and \eqref{MWstrong} are a weakening of the first well known condition by Gordin \cite{G69}:  $\sum_{k \geq 1} \|{\mathbb E}(X_k | {\mathcal F}_0)\|_2 < \infty$, under which a martingale + coboundary decomposition was established. Note that, under \eqref{MWstrong} the unique variable $D$ satisfying ${\mathbb L}^2_0$-MMA is given by 
\begin{equation}\label{Gmart}
D= \sum_{k=0}^\infty P_0(X_k) \, , \quad \text{where $P_0(X_k)={\mathbb E}(X_k|{\mathcal F}_0)-{\mathbb E}(X_k|{\mathcal F}_{-1}),$}
\end{equation}
which is  exactly the martingale given in Gordin's paper \cite{G69}. Note also that, under \eqref{MWstrong}, the series \eqref{Gmart} is normally convergent in ${\mathbb L}^2$.

To conclude this introduction, note that all the results presented in this note extend naturally to variables with values in a Hilbert space (by replacing absolute values and products of variables by the Hilbertian norm and the scalar product in the appropriate places), since the proofs rely on orthogonality arguments in ${\mathbb L}^2$.

\section{The condition \eqref{GL} and the ${\mathbb L}^2$-MA property}
\setcounter{equation}{0}
Let $P_0$ be the operator defined in \eqref{Gmart}. Let us first recall a lemma stated in \cite{D15}.
\begin{Lemma}\label{D15}
Assume that \eqref{GL} holds and let $D_n=P_0(X_0+S_n)$. Then 
$$
\lim_{n \rightarrow \infty} \frac{\|{\mathbb E}(S_n|{\mathcal F}_0)\|_2}{\sqrt n}=0 \quad \text{and} \quad 
\lim_{n \rightarrow \infty} \frac 1 n \left \| S_n - \sum_{k=1}^n D_n\circ T^k \right \|_2^2 = 0 \, .
$$
\end{Lemma}
From Lemma \ref{D15}, we see that the asymptotic distribution of $n^{-1/2} S_n$ is equivalent to that of $n^{-1/2} \sum_{k=1}^n D_n\circ T^k$. Lemma \ref{D15}  (and its proof) is the main tool to get the following Proposition.
\begin{Proposition} \label{mainProp}
Assume that \eqref{GL} holds. The following statements are equivalent:
\begin{enumerate}
\item The ${\mathbb L}^2$-\textup{MA} property holds.
\item The series \eqref{Gmart} converges in ${\mathbb L}^2$.
\item $ \displaystyle  \lim_{N \rightarrow \infty} \sup_{n \geq N} \frac 1 n \sum_{k=1}^n {\mathbb E}(X_0{\mathbb E}(S_{k-1}|{\mathcal F}_{-N})) =0 \, .$
\end{enumerate}
Moreover, if ${\mathbb L}^2$-\textup{MA} holds, the approximating martingale is Gordin's martingale \eqref{Gmart}.
\end{Proposition}
As an immediate consequence of Proposition \ref{mainProp}, we get the following corollary:
\begin{corollary}\label{mainCor}
If \eqref{GL} holds and  $\displaystyle \lim_{N, n \rightarrow \infty}  | {\mathbb E}(X_0{\mathbb E}(S_{n}|{\mathcal F}_{-N}))| =0 $, then  ${\mathbb L}^2$-\textup{MA} holds with $D$ defined by \eqref{Gmart}.
\end{corollary}

\noindent {\bf Example 1: $\alpha$-dependent sequences.} Let $F(t)={\mathbb P}(X_0 \leq t)$ and let $Q$ be the cadlag inverse of $H=1-F$. Define the $\alpha$-dependent coefficients of the sequence $(X_i)_{i \in {\mathbb Z}}$ as in \cite{R00} (see equation (2.30) there):
$$
\alpha(k)=\sup_{t \in {\mathbb R}} \|{\mathbb E}({\bf 1}_{X_k \leq t}|{\mathcal F}_0)-F(t)\|_1 \, .
$$
Proceeding exactly as in the proof of Theorem 4.2 in \cite{R00}, we see that both  $\displaystyle \lim_{N, n \rightarrow \infty} | {\mathbb E}(X_0{\mathbb E}(S_{n}|{\mathcal F}_{-N}))| =0 $  and \eqref{GL} are satisfied as soon as 
\begin{equation}\label{Rio}
\sum_{k \geq 0} \int_0^{\alpha(k)} Q^2(u) du < \infty \, .
\end{equation}
Hence, by Corollary \ref{mainCor}, we obtain that ${\mathbb L}^2$-MA holds under \eqref{Rio}. Note that condition \eqref{Rio} is Rio's condition for the CLT of $\alpha$-dependent sequence (see Theorem 4.2 in \cite{R00}, which is stated in the ergodic case). As recalled in the introduction, the property ${\mathbb L}^2$-MA implies  that $n^{-1/2} S_n $ converges in distribution to $\sqrt{{\mathbb E}(D^2|{\mathcal I})} Z$, where $Z$  is ${\mathcal N}(0,1)$-distributed and independent of ${\mathcal I}$, so we recover Rio's CLT \cite{R00}. Identifying the variances in the CLT, one can also prove that, under \eqref{Rio},
$$
{\mathbb E}(D^2|{\mathcal I})= \sum_{k \in {\mathbb Z}} {\mathbb E}(X_0X_k|{\mathcal I}) \ \  {a.s.}
$$
To conclude, note that the quenched WIP under \eqref{Rio} has been proved in \cite{DMP14}.

\medskip

\noindent{\bf Proof of Proposition \ref{mainProp}.} 
From Lemma \ref{D15}, the  ${\mathbb L}^2$-MA property holds if and only if there exists $D$ in ${\mathbb L}^2$ such that ${\mathbb E}(D \circ T|{\mathcal F_0})= 0 $ a.s., and
\begin{equation}\label{convDn}
\lim_{n \rightarrow} \frac 1 n \left \|  \sum_{k=1}^n (D_n- D)\circ T^k \right \|_2^2 = 0 \, .
\end{equation}
Now, by orthogonality, 
$$
\frac 1 n \left \|  \sum_{k=1}^n (D_n- D)\circ T^k \right \|_2^2=\|D_n-D\|_2^2 \, , 
$$
so that \eqref{convDn} holds if and only if $D_n$ converges in ${\mathbb L}^2$. This proves the equivalence of 1. and 2. 

The main point is to prove that 2. is equivalent to 3.  By the Cauchy criterion, 2. is equivalent to 
\begin{equation}\label{Cauchy}
\lim_{N\rightarrow \infty} \sup_{n\geq N} \left \| \sum_{k=N}^n P_0(X_k) \right \|_2^2 =0 \, .
\end{equation}
For any positive integer $N$, let $X_{k, N}={\mathbb E}(X_k|{\mathcal F}_{k-N})$, so that $(X_{n, N})_{k\in {\mathbb Z}}$ is a strictly stationary sequence. It is then easy to see that 
$$
P_0(X_{k, N}) = P_0(X_k) \quad \text{if $k \geq N$, and} \quad P_0(X_{k, N}) =0 \ \text{a.s, if $k < N$.}
$$
It follows that 
$$
\sum_{k=N}^n P_0(X_k) = P_0(S_{n, N})\, , \quad \text{where} \quad S_{n,N}=X_{1,N}+ \cdots + X_{n,N} \, .
$$
Let then $D_{n, N}=P_0(X_{0,N}+S_{n, N})$. From the proof of Lemma \ref{D15} (see the proof of Lemma 3.1 in \cite{D15}), we have 
\begin{equation}\label{mainterm}
\frac 1 n \left \| S_{n,N} - \sum_{k=1}^n D_{n,N}\circ T^k \right \|_2^2 \leq \sum_{i=1}^n \|{\mathbb E}(X_{i,N}|{\mathcal F}_0)\|_2^2 + \frac 1 n  \sum_{i=1}^n  \left ( \frac 1 i \sum_{p=1}^i \beta_{p,N}  \right ) \, ,
\end{equation}
where $\beta_{p,N}= \sum_{i=p}^\infty i \|P_0(X_{i, N})\|_2^2$. Here note that $\beta_{p,N} \leq \beta_p= \sum_{i=p}^\infty i \|P_0(X_{i})\|_2^2$ (which is finite because \eqref{GL} is equivalent to ${\mathbb E}(X_0|{\mathcal F}_{-\infty})=0$ and $\sum_{i=1}^\infty i \|P_0(X_{i})\|_2^2 < \infty$; see Remark \ref{Rq2} of Section \ref{Sec4} for more details). Since $\beta_p \rightarrow 0$ as $p \rightarrow \infty$, we deduce that 
\begin{equation}\label{step1}
\lim_{N\rightarrow \infty} \sup_{n\geq N} \frac 1 n  \sum_{i=1}^n  \left ( \frac 1 i \sum_{p=1}^i \beta_{p,N}  \right ) =0 \, .
\end{equation}
Now, we also have that $ \|{\mathbb E}(X_{i,N}|{\mathcal F}_0)\|_2^2 \leq \|{\mathbb E}(X_{i}|{\mathcal F}_0)\|_2^2 $, and $\lim_{N \rightarrow \infty} \|{\mathbb E}(X_{i,N}|{\mathcal F}_0)\|_2^2= 0$. Hence, by the dominated convergence theorem, since \eqref{GL} holds, we have 
\begin{equation}\label{step2}
\lim_{N\rightarrow \infty} \sum_{i=1}^\infty \|{\mathbb E}(X_{i,N}|{\mathcal F}_0)\|_2^2 =0 \, . 
\end{equation}
From \eqref{mainterm}, \eqref{step1} and \eqref{step2}, we infer that 
\begin{equation}\label{mainlim}
\lim_{N\rightarrow \infty} \sup_{n\geq N}
\frac 1 n \left \| S_{n,N} - \sum_{k=1}^n D_{n,N}\circ T^k \right \|_2^2 =0 \, .
\end{equation}
Since 
$$
 \frac 1 n \left \| \sum_{k=1}^n D_{n,N} \circ T^k \right \|_2^2= \|D_{n,N}\|_2^2= \left  \| \sum_{k=N}^n P_0(X_k) \right \|_2^2 \, ,
$$
we infer from \eqref{mainlim} that
\begin{equation}\label{mainlim2}
\lim_{N\rightarrow \infty} \sup_{n\geq N}
  \left | \frac{\|S_{n,N}\|_2^2}{n} - \left  \| \sum_{k=N}^n P_0(X_k) \right \|_2^2 \right | =0 \, .
\end{equation}
From \eqref{Cauchy} and \eqref{mainlim2}, we see that 2. is in fact equivalent to 
\begin{equation}\label{Cauchy2}
\lim_{N\rightarrow \infty} \sup_{n\geq N} \frac{\|S_{n,N}\|_2^2}{n} = 0 \, .
\end{equation}
Now, it is an easy exercise to show that
$$
\frac{\|S_{n,N}\|_2^2}{n}=\gamma_N(0) + \frac 2 n \sum_{k=1}^n \sum_{i=1}^{k-1} \gamma_N(i) \, , \quad \text{where, for $i\geq 0$,} \quad 
\gamma_N(i)={\mathbb E}(X_{0,N}X_{i,N})={\mathbb E}(X_0{\mathbb E}(X_{i}|{\mathcal{F}_{-N}})) \, .
$$
Since $\gamma_N(1)+ \cdots + \gamma_N(k-1)= {\mathbb E}(X_0{\mathbb E}(S_{k-1}|{\mathcal{F}_{-N}}))$, we get that 
\begin{equation}\label{endterm}
\frac{\|S_{n,N}\|_2^2}{n}=\|{\mathbb E}(X_{0}|{\mathcal{F}_{-N}}) \|_2^2 + \frac 2 n \sum_{k=1}^n {\mathbb E}(X_0{\mathbb E}(S_{k-1}|{\mathcal{F}_{-N}})) \, .
\end{equation}
Since \eqref{GL} holds, we have $\lim_{N\rightarrow \infty}\|{\mathbb E}(X_{0}|{\mathcal{F}_{-N}}) \|_2^2=0$, and we infer from \eqref{endterm} that \eqref{Cauchy2} is equivalent to 3. Since we have proved that \eqref{Cauchy2} is also equivalent to 2., we conclude that 3. is equivalent to 2., and the proof of Proposition \ref{mainProp} is complete. $\square $

\section{The condition \eqref{GL} and the ${\mathbb L}^2_0$-MA property}\label{L20MA}

\setcounter{equation}{0}

In this section and the next one, we shall need the following lemma:
\begin{Lemma}\label{suiteLemme}
Assume that \eqref{GL} holds. Then
$$
\lim_{n \rightarrow \infty}\frac{1}{\sqrt n} \left \| \max_{1 \leq k \leq n} | {\mathbb E}(S_k|{\mathcal F}_0)|\right \|_2= 0 \quad \text{and} \quad \lim_{n \rightarrow \infty}\frac{1}{\sqrt n}  \max_{1 \leq k \leq n} | {\mathbb E}(S_k|{\mathcal F}_0)|= 0 \ \ a. s.
$$
\end{Lemma}

From Lemma \ref{suiteLemme} and Proposition 7 of Cuny and Peligrad \cite{CP12}, we get the following proposition.

\begin{Proposition}\label{Prop2}
Assume that \eqref{GL} holds. If the series \eqref{Gmart} converges in ${\mathbb L}^2$ and
\begin{equation}\label{CPcond}
\lim_{N \rightarrow \infty} {\mathbb E} \left ( \sup_{k \geq N} \left | \sum_{i=k}^\infty P_0(X_i) \right |^2\right )= 0 \, ,
\end{equation}
then the ${\mathbb L}^2_0$-\textup{MA} property holds. Moreover, the approximating martingale is Gordin's martingale \eqref{Gmart}.
\end{Proposition}

\begin{Remark} By the monotone convergence theorem, 
is is easy to see that condition \eqref{CPcond} implies that the series \eqref{Gmart} converges almost surely. We shall see in Section \ref{QLin} that \eqref{CPcond} is not necessary for ${\mathbb L}^2_0$-\textup{MA} to hold.
\end{Remark}

\noindent {\bf Example 2: causal linear processes with martingale difference innovations.} For $k \in {\mathbb Z}$, let $\varepsilon_{k}= \varepsilon_0 \circ T^k$, where $\varepsilon_0$ is square integrable, ${\mathcal F}_0$ -measurable, and such that ${\mathbb E}(\varepsilon_0|{\mathcal{F}}_{-1})=0$ a.s. and Var$(\varepsilon_0)>0$. Let also $(a_k)_{k \geq 0}$ be a sequence in $\ell^2$. Define then  the stationary sequence 
$$
X_k=X_0 \circ T^k = \sum_{i=0}^{\infty} a_i \varepsilon_{k-i} \, .
$$
One can easily check that $P_0(X_k)=a_k \varepsilon_0$ and that the condition \eqref{GL} is exactly $\sum_{i \geq 0} i a_i^2 < \infty$ (see Remark \ref{Rq2} below for more explanations about this equivalence). If $\sum_{i \geq 0} i a_i^2 < \infty$, the four following points are equivalent:
\begin{enumerate}
\item $n^{-1/2}S_n$ converges in distribution.
\item The ${\mathbb L}^2$-MA property holds.
\item The ${\mathbb L}^2_0$-MA property holds.
\item The series $\sum_{k=0}^\infty a_k $ converges.
\end{enumerate}
Moreover, if ${\mathbb L}^2$-MA holds, the approximating martingale is $D=(\sum_{k \geq 0} a_k)\varepsilon_0$. 

The equivalence of items 2.-3.-4. follows easily from Proposition \ref{mainProp} and Proposition \ref{Prop2}. Let us briefly explain why item 1. is equivalent to item 4. As we already mentioned, it follows from Lemma \ref{D15} that the asymptotic distribution of $n^{-1/2}S_n$ is equivalent to that of 
$(\sum_{k=0}^n a_k) n^{-1/2}S_n(\varepsilon)$, where $S_n(\varepsilon)= \varepsilon_1 + \cdots + \varepsilon_n$. From the CLT for sums of stationary martingale differences, we  deduce that the asymptotic distribution of $n^{-1/2}S_n$ is equivalent to that of $(\sum_{k=0}^n a_k) \sqrt{{\mathbb E}(\varepsilon_0^2|{\mathcal I})} Z $, where $Z$ is ${\mathcal N}(0,1)$-distributed and independent of ${\mathcal I}$.  Since $a_n \rightarrow 0$ as $n \rightarrow \infty$, we easily deduced that $n^{-1/2}S_n$ converges in distribution iff $\sum_{k=0}^\infty a_k $ converges.

As a concluding remark concerning this example, let us mention that Voln\'y and Woodroofe \cite{VW10} constructed an example of a causal linear process with iid centered and square integrable innovations, for which $a_k\geq 0$ for all $k\geq0$,  $ 0<\sum_{k\geq 0} a_k <\infty $, and the sequence $n^{-1/2}S_n$ satisfies the CLT (and even the ${\mathbb L}^2$-MMA property) but not the quenched CLT (hence ${\mathbb L}^2_0$-MA does not hold).

\medskip

\noindent{\bf Proof of Lemma \ref{suiteLemme}.} Clearly, it suffices to prove that 
\begin{equation}\label{step1bis}
\lim_{m \rightarrow \infty} \limsup_{n \rightarrow \infty} \frac 1 n \left \| \max_{m \leq k \leq n} | {\mathbb E}(S_k-S_m|{\mathcal F}_0)|\right \|_2^2= 0 \, ,
\end{equation}
and 
\begin{equation}\label{step2bis}
\lim_{m \rightarrow \infty} \limsup_{n \rightarrow \infty} \frac 1 n  \max_{m \leq k \leq n} | {\mathbb E}(S_k-S_m|{\mathcal F}_0)|^2= 0 \quad \text{almost surely.}
\end{equation}
Now, for $m \leq k \leq n $,
$$({\mathbb E}(S_k-S_m|{\mathcal F}_0))^2 = \sum_{i=m+1}^k \sum_{j=m+1}^k {\mathbb E}(X_i|{\mathcal F}_0){\mathbb E}(X_j|{\mathcal F}_0)
\leq k \sum_{i=m+1}^k 
({\mathbb E}(X_i|{\mathcal F}_0))^2 \leq n \sum_{i=m+1}^\infty 
({\mathbb E}(X_i|{\mathcal F}_0))^2 \, ,
$$
where the first upper bound  holds because 
$$|{\mathbb E}(X_i|{\mathcal F}_0){\mathbb E}(X_j|{\mathcal F}_0)|
\leq \frac 1 2 ({\mathbb E}(X_i|{\mathcal F}_0))_2^2+ \frac 1 2({\mathbb E}(X_j|{\mathcal F}_0))_2^2 \, .
$$ 
Hence, both \eqref{step1bis} and \eqref{step2bis} follow easily from  the fact that
$\sum_{n \geq 0} \|{\mathbb E}(X_n|{\mathcal F}_0)\|_2^2 < \infty$. $\square$

\medskip

\noindent{\bf Proof of Proposition \ref{Prop2}.} By Lemma \ref{suiteLemme}, it suffices to prove that there exists a  ${\mathcal F}_0$-measurable random variable $D$ in ${\mathbb L}^2$ such that ${\mathbb E}(D \circ T|{\mathcal F_0})= 0 $ a.s. and 
$$
\lim_{n \rightarrow \infty} \frac 1 n {\mathbb E} \left ( \left ( S_n-{\mathbb E}(S_n|{\mathcal F}_0) - \sum_{k=1}^n D\circ T^k \right )^2 \Big | {\mathcal F}_0 \right ) = 0 \, , \text{almost surely and in ${\mathbb L}^1$.}
$$
The ${\mathbb L}^1$ convergence follows from Proposition \ref{mainProp}. To prove the almost sure convergence, we can proceed exactly as for the proof of Proposition 7 in \cite{CP12}. In fact, in the context of Markov chain, our condition \eqref{CPcond} is exactly the condition  (14) of Cuny and Peligrad \cite{CP12}. It is easy to see that their proof still works in the general context of strictly stationary sequences. $\square$

\section{The condition \eqref{GL} and the ${\mathbb L}^2_0$-MMA property}\label{Sec4}

\setcounter{equation}{0}

The sequence $(X_i, {\mathcal F}_i)_{i \geq 0}$ is said to satisfy Hannan's condition if
\[
{\mathbb E}(X_0 | {\mathcal F}_{- \infty})= 0 \ \text{a.s. and} \quad
\sum_{k \geq 0} \|P_0(X_k)\|_2 < \infty \, .  \tag{H}  \label{H}
\]
This condition was first introduced by Hannan \cite{H73}, to prove the CLT for $n^{-1/2}S_n $ and some linear statistics based on $X_1, \ldots, X_n$ (assuming an extra condition stronger than ergodicity). In 1979, Hannan \cite{H79} proved that the WIP holds under the same condition. We refer to the paper by Dedecker et al. \cite{DMV07} for more comments about condition (H). In 2011, Gordin and Peligrad \cite{GP11} showed that ${\mathbb L}^2$-MMA holds under \eqref{H} for $D$ defined by \eqref{Gmart}. The fact that, under \eqref{H}, ${\mathbb E}(D^2|{\mathcal I})=\sum_{k \in {\mathbb Z}} {\mathbb E}(X_0X_k|{\mathcal I}) $ is proved in Corollary 1 of \cite{DMV07}. 

Concerning the ${\mathbb L}^2_0$-MMA property, the following result has been proved by Cuny and Voln\'y \cite{CV13}: if \eqref{H} is satisfied, then
\begin{equation}\label{partialR}
\lim_{n \rightarrow \infty} \frac 1 n {\mathbb E} \left (  \max_{1 \leq k \leq n} \left( S_k-{\mathbb E}(S_k| {\mathcal F}_0) - \sum_{i=1}^k D\circ T^i \right )^2 \Big | {\mathcal F}_0 \right ) = 0 \, , \text{almost surely.}
\end{equation}

Combining the result by Gordin and Peligrad \cite{GP11}, the approximation \eqref{partialR} and Lemma \ref{suiteLemme}, we obtain the following result:

\begin{Proposition}\label{maxProp}
If \eqref{GL} $+$ \eqref{H} holds, then the ${\mathbb L}^2_0$-\textup{MMA} property holds. Moreover, the approximating martingale is Gordin's martingale \eqref{Gmart}.
\end{Proposition}

\begin{Remark}\label{Rq2} Let $P_i(X_k)={\mathbb E}(X_k|{\mathcal F}_i)-{\mathbb E}(X_k|{\mathcal F}_{i-1})$.
If ${\mathbb E}(X_0|{\mathcal F}_{-\infty})=0$ almost surely, then ${\mathbb E}(X_k|{\mathcal F}_0)= \sum_{i \leq 0} P_i(X_k)$. Hence, by orthogonality and stationarity, we infer that
$\|{\mathbb E}(X_k|{\mathcal F}_0)\|_2^2=\sum_{i \geq k} \|P_0(X_i)\|_2^2$. It follows that 
the condition \eqref{GL} $+$ \eqref{H} may be written with the help of the projection operator $P_0$ only, as follows: 
\begin{equation}
{\mathbb E}(X_0 | {\mathcal F}_{- \infty})= 0 \ \text{a.s.}\, , \quad \sum_{k \geq 0} k\|P_0(X_k)\|_2^2 < \infty \quad  \text{and} \quad
\sum_{k \geq 0} \|P_0(X_k)\|_2 < \infty \, .    \label{GLH}
\end{equation}
As quoted in \cite{DMV07}, a condition expressed in terms of $\|P_0(X_k)\|_2$ is particularly interesting for applications to functions of linear processes with iid innovations. We shall apply it to a larger class of examples in Section \ref{QLin}.
\end{Remark}

As a consequence of Proposition \ref{maxProp} and Remark \ref{Rq2}, we obtain the following corollary:

\begin{corollary}\label{dec}
Assume that ${\mathbb E}(X_0 | {\mathcal F}_{- \infty})= 0$ a.s. and that $\|P_0(X_i)\|_2$ is bounded by $u_i$, where $(u_i)_{i \geq 0}$ is nonincreasing and summable. Then the ${\mathbb L}^2_0$-\textup{MMA} property holds.
\end{corollary}

\noindent{\bf Proof of Corollary \ref{dec}.} It suffices to notice that if $(u_i)_{i \geq 0}$ is nonincreasing and summable, then $\sum_{k\geq 0} k u_k^2 < \infty$, and to apply Proposition \ref{maxProp} and Remark \ref{Rq2}. $\square$.

\medskip 

Recall that the condition \eqref{MWstrong} of the introduction (that is  $ \sum_{k>0} k^{-1/2} \|{\mathbb E}(X_k|{\mathcal F}_0)\|_2 < \infty$) also implies ${\mathbb L}^2_0$-MMA. In the next proposition, we give the hierarchy between \eqref{MWstrong} and \eqref{GL} + \eqref{H}.

\begin{Proposition}\label{comp}
We have the implication: \eqref{MWstrong} $\Rightarrow$ \eqref{GL} $+$ \eqref{H}. Moreover:
\begin{enumerate}
\item There exists a causal linear process with iid innovations satisfying \eqref{GL} $+$ \eqref{H} but not \eqref{MWstrong}.
\item  If the sequence $(\|P_0(X_i)\|_2)_{i \geq 0}$ is nonincreasing then the three conditions \eqref{H}, \eqref{GL} $+$ \eqref{H} and \eqref{MWstrong} are equivalent.
\end{enumerate}
\end{Proposition}

To prove Proposition \ref{comp}, we will need the following lemma, which has interest in itself. 

\begin{Lemma}\label{Lemmeseries}
Let $(u_i)_{i\geq 0}$ be a sequence of nonnegative numbers, and  $q \in (1, \infty)$. 
Then
\begin{equation}\label{series}
 \sum_{k>0} \left ( \frac 1 k   \sum_{i=k}^\infty u_i^q \right)^{1/q}  < \infty \quad \Longrightarrow  \quad \sum_{k>0} u_k < \infty \, .
\end{equation}
Moreover, if $(u_i)_{i\geq 0}$ is nonincreasing, the two conditions are equivalent. 
\end{Lemma}

\noindent{\bf Proof of Lemma \ref{Lemmeseries}.} The property \eqref{series} is stated in Lemma A1 of \cite{DMP11} (in fact, a more general result is proved there). It remains to prove the converse implication, in the case where $(u_i)_{i\geq 0}$ is nonincreasing. 

Let $v_j=u_{2^j}$. In the case where $(u_i)_{i\geq 0}$ is nonincreasing, it is easy to see that the summability of $(u_i)_{i\geq 0}$ is equivalent to the summability of $(2^j v_j)_{j \geq 0}$. Let $k \in \{2^{l}, \ldots, 2^{l+1}-1\}$. Then 
$1/k \leq 2^{-l}$ and
$$
\sum_{i=k}^\infty u_i^q \leq \sum_{i=2^l}^\infty u_i^q \leq \sum_{j=l}^\infty 2^j v_j^q \, .
$$
It follows that 
$$
\sum_{k>0} \left ( \frac 1 k   \sum_{i=k}^\infty u_i^q \right)^{1/q} \leq \sum_{l=0}^\infty 2^l \left ( \frac 1 {2^l}   \sum_{j=l}^\infty 2^j v_j^q \right)^{1/q} \leq 
\sum_{l=0}^\infty 2^{l(1-1/q)}    \sum_{j=l}^\infty 2^{j/q} v_j \, .
$$
Interverting the sums in the last bound, and noting that 
$\sum_{l=0}^j 2^{l(1-1/q)} \leq C 2^{j(1-1/q)}$ for some $C>0$, we get 
$$
\sum_{k>0} \left ( \frac 1 k   \sum_{i=k}^\infty u_i^q \right)^{1/q} \leq C \sum_{j=l}^\infty 2^{j} v_j < \infty \, .
$$
This completes the proof. $\square$

\medskip

\noindent{\bf Proof of Proposition \ref{comp}.} The fact that \eqref{MWstrong} $\Rightarrow$  \eqref{H} is well known; it follows from Remark \ref{Rq2} and \eqref{series} of Lemma \ref{Lemmeseries} by taking  $q=2$ and
$u_i=\|P_0(X_i)\|_2$. The fact that \eqref{MWstrong} $\Rightarrow$  \eqref{GL} (already mentioned in the introduction) follows by noting that $(k^{-1/2} \|{\mathbb E}(X_k|{\mathcal F}_0)\|_2)_{k >0}$ is nonincreasing, and by using the argument of the proof of Corollary \ref{dec}. The point 2. of Proposition \ref{comp} follows from the last assertion of Lemma \ref{Lemmeseries} by taking again 
$u_i=\|P_0(X_i)\|_2$ (so that $(u_i)_{i\geq 0}$ is assumed to be nonincreasing). 

It remains to prove 1. Let 
$$
X_k= \sum_{i=0}^{\infty} a_i \varepsilon_{k-i} \, ,
$$
where $(\varepsilon_i)_{i \in {\mathbb Z}}$ is a sequence of iid centered and square integrable random variables such that Var$(\varepsilon_0)>0$, and $(a_i)_{i \geq 0}$ belongs to $\ell^2$. Let ${\mathcal F}_k= \sigma(\varepsilon_i, i \leq k)$. As in Example 2 of Section \ref{L20MA}, we see that $P_0(X_k)=a_k \varepsilon_0$. Define now the sequence $(a_i)_{i \geq 0}$ as follows: 
$$
a_i=0  \ \text{if $i\neq 2^k$, and} \ a_{2^k}= \frac{1}{2^{k/2} k^b  } \, , \ \text{for $b \in (1/2,1)$. }
$$
From the definition of the $a_i$'s, we immediately see that 
$\sum_{i\geq 0} a_i < \infty$ and $\sum_{i\geq 0} ia_i^2 < \infty$, so that both \eqref{H} and \eqref{GL} hold. To prove that \eqref{MWstrong} does not hold, it suffices to prove that
$$
  \sum_{k>0} \frac{1}{\sqrt k} \left (\sum_{i=k}^{\infty} a_i^2 \right )^{1/2} = \infty \, .
$$
Write then 
\begin{equation}\label{point1}
\sum_{k>1} \frac{1}{\sqrt k} \left (\sum_{i=k}^{\infty} a_i^2 \right )^{1/2}= \sum_{l=1}^\infty \sum_{k=2^{l-1}+1}^{2^{l}} \frac{1}{\sqrt k} \left (\sum_{i=k}^{\infty} a_i^2 \right )^{1/2} \, .
\end{equation}
Now, for $k \in \{2^{l-1}+1, \ldots, 2^{l}\}$, 
\begin{equation}\label{point2}
    \left (\sum_{i=k}^{\infty} a_i^2 \right )^{1/2} \geq a_{2^l} \, .
\end{equation}
Moreover, there exists $c>0$ such that 
\begin{equation}\label{point3}
 \sum_{k=2^{l-1}+1}^{2^{l}} \frac{1}{\sqrt k} \geq c 2^{l/2} \, .
\end{equation}
From \eqref{point1}, \eqref{point2} and \eqref{point3}, we see that 
$$
\sum_{k>1} \frac{1}{\sqrt k} \left (\sum_{i=k}^{\infty} a_i^2 \right )^{1/2} \geq \sum_{l=1}^\infty \frac{c}{l^b}= \infty \, ,
$$
which proves that \eqref{MWstrong} does not hold. $\square$

\section{Semi-linear processes and H\"older functions of semi-linear processes}\label{QLin}

\setcounter{equation}{0}

In this section, we revisit and develop a class of examples originally introduced in \cite{D15} to show that the \eqref{GL} condition is not sufficient to guarantee the CLT (see Proposition \ref{CE} below).

Let $({\mathcal X}, {\mathcal B}, \mu)$
be a probability space. Let $(\alpha_i)_{i \geq 0}$
be a sequence of functions in ${\mathbb L}^2({\mathcal X})$   such that
\begin{equation}\label{alpha}
\mu(\alpha_i)=0  \quad  \text{and} \quad
 \sum_{j \geq 0} \|\alpha_j\|_{2, {\mathcal X}}^2 < \infty\, .
\end{equation}

We consider the  space $\Omega={\mathcal X}^{\otimes {\mathbb Z}}$ and the probabilty
${\mathbb P}= \mu^{\otimes {\mathbb Z}}$. 
The transformation $T$ is the shift on $\Omega$ defined by
$(T(\omega))_i= \omega_{i+1}$. Clearly ${\mathbb P}$ is invariant by $T$ and the couple
$(T, {\mathbb P})$ is ergodic.

Starting from the sequence $(\alpha_i)_{i \geq 0}$ and from the projections $\pi_j(\omega)=\omega_j$, we define the sequence
$(A_i)_{i\geq 0}$ of functions of ${\mathbb L}^2$ by:
$A_i= \alpha_i \circ \pi_0$. The sequence $(X_i)_{i \in {\mathbb Z}}$
is then defined by:
\begin{equation}\label{QLdef}
X_0 =\sum_{j=0}^\infty A_j \circ T^{-j}= \sum_{j=0}^\infty \alpha_j \circ \pi_{-j} \quad \mbox{and} \quad X_i= X_0 \circ T^i =\sum_{j=0}^\infty A_j \circ T^{i-j} \, .
\end{equation}
Note that these series are well defined in ${\mathbb L}^2({\mathbb P})$ because
$(A_j \circ T^{i-j})_{j \geq 0}$ is a sequence of independent and centered random variables, and
$$\sum_{j\geq 0} \|A_j \circ T^{i-j}\|_2^2 = \sum_{j\geq 0} \|A_j \|_2^2= \sum_{j \geq 0} \|\alpha_j\|_{2, {\mathcal X}}^2 < \infty \, . $$
We say that the stationary and ergodic sequence $(X_i)_{i \in {\mathbb Z}}$ defined by \eqref{QLdef} is a {\em semi-linear process}. Note that the class of semi-linear processes contains the class of causal linear processes with iid innovations (it suffices to consider the particular case where $\alpha_i=a_i \varepsilon_0$, with 
$\varepsilon_0 \in {\mathbb L}^2({\mathcal X})$, $\mu(\varepsilon_0)=0$, and $(a_i)_{i\geq 0}\in \ell^2$).

Let ${\mathcal F}_i= \sigma( \pi_j, j \leq i)$.
Clearly, $X_0$ is ${\mathcal F}_0$-measurable and
${\mathcal F}_{-\infty}$ is ${\mathbb P}$-trivial by the $0-1$ law. We see that
$P_0(X_i)=A_i$ for any positive integer $i$. By Remark \ref{Rq2}, we see that \eqref{GL} holds iff 
\begin{equation}\label{GL2}
 \sum_{j\geq 0} j \|A_j \|_2^2=\sum_{j \geq 0}  j \|\alpha_j\|_{2, {\mathcal X}}^2 < \infty \, .
\end{equation}

\subsection{The ${\mathbb L}^2_0$-MA property for semi-linear processes}\label{L20MAbis}

For the class of semi-linear processes, we obtain a more precise result than Proposition \ref{Prop2}.

\begin{Proposition} \label{mainProp2} Let $(X_i)_{i \in {\mathbb Z}}$ be a semi-linear process, as in \eqref{QLdef}.
Assume that \eqref{GL} holds (i.e. \eqref{GL2} holds). The following statements are equivalent:
\begin{enumerate}
\item The ${\mathbb L}^2$-\textup{MA} property holds.
\item The ${\mathbb L}^2_0$-\textup{MA} property holds.
\item The series $\sum_{i\geq 0} A_i$ converges in ${\mathbb L}^2$.
\end{enumerate}
Moreover, if ${\mathbb L}^2$-\textup{MA} holds, the approximating martingale is Gordin's martingale $D=\sum_{i\geq 0} A_i$.
\end{Proposition}

\begin{Remark}
Let $(X_i)_{i \in {\mathbb Z}}$ be a semi-linear process, as in \eqref{QLdef}.
Assume that \eqref{GL} holds, and let $D_n=\sum_{i= 0}^n A_i$. Using the same arguments as in the proof of Proposition \ref{mainProp2}, we can prove that $n^{-1/2}S_n$ satisfies the Quenched CLT with convergence to a ${\mathcal N}(0, \sigma^2)$ distribution
iff $n^{-1/2} \sum_{i=1}^n D_n\circ T^i$ converges in distribution to ${\mathcal N}(0, \sigma^2)$. Since the variables $(D_n \circ T^i)_{1 \leq i \leq n}$ are independent, we infer from Lindeberg's theorem that $n^{-1/2}S_n$ satisfies the Quenched CLT as soon as ${\mathbb E}(D_n^2)$ converges to $\sigma^2$ and, for any $\epsilon >0$,
$\lim_{n \rightarrow \infty}{\mathbb E}(D_n^2{\bf 1}_{|D_n|\geq \sqrt n \epsilon})=0$. Of course, these conditions are satisfied as soon as $D_n$ converges to $D$ in ${\mathbb L}^2$, with $\sigma^2={\mathbb E}(D^2)$.
\end{Remark}

Concerning the optimality of Proposition \ref{mainProp2}, let us recall the counter-examples given in \cite{D15}. 

\begin{Proposition}\label{CE}
Let $(v_n)_{n \geq 0}$ and $(u_n)_{n \geq 0}$ be two sequences of nonnegative numbers such that
$$
v_n \geq n \, , \quad \sum_{n=1}^\infty \frac 1{v_n} = \infty \, ,\quad  \sum_{n=0}^\infty u_n = \infty \, , \quad  \text{and} \quad \sum_{n=0}^\infty v_n u_n^2 < \infty \, .
$$
\begin{enumerate}
\item There exists a causal linear process $(X_i)_{i \in {\mathbb Z}}$ with iid innovations such that $\|P_0(X_i)\|_2=u_i$, and $n^{-1/2} S_n$ is not stochastically bounded.
\item There exists a semi-linear process $(X_i)_{i \in {\mathbb Z}}$ such that $\|P_0(X_i)\|_2 \leq u_i$ and $n^{-1}{\mathbb E}(S_n^2) \rightarrow 1$ as $n\rightarrow \infty$, but $n^{-1/2}S_n$ does not converge in distribution. 
\end{enumerate}
\end{Proposition}

\begin{Remark}
Taking $v_n=n \log (n+e)$ (so that $\sum_{n>0} v_n^{-1}= \infty$), we  infer from Proposition \ref{CE} that the condition ${\mathbb E}(X_0 | {\mathcal F}_{- \infty})= 0 \ \text{a.s.}$ $+$
$\sum_{i >0} i \log(i)  \|P_0(X_i)\|_2^2 < \infty$ (which is equivalent to \eqref{GLlog})
is not sufficient for the CLT (with normalization $\sqrt n$), even if $(X_i)_{i \in {\mathbb Z}}$ is a causal linear process with iid innovations.
\end{Remark}

\noindent{\bf Proof of Proposition \ref{mainProp2}.}
By Proposition \ref{mainProp}, 1. and 3. are equivalent under \eqref{GL}, and the approximating martingale is $D=\sum_{i\geq 0} A_i$. It remains to prove that if 1. holds, then 2.  also holds with  the same $D$.

Since \eqref{GL} holds,  by Lemma \ref{suiteLemme}, Item 2. with $D=\sum_{i\geq 0} A_i$  is equivalent to 
\begin{equation}\label{partialR2}
\lim_{n \rightarrow \infty} \frac 1 n {\mathbb E} \left (   \left( S_n-{\mathbb E}(S_n| {\mathcal F}_0) - \sum_{i=1}^n D\circ T^i \right )^2 \Big | {\mathcal F}_0 \right ) = 0 \, , \, \text{almost surely and in ${\mathbb L}^1$.}
\end{equation}
Note that the variables $S_n-{\mathbb E}(S_n| {\mathcal F}_0)$ and $D\circ T^i$ for $i\geq 1$ are independent of ${\mathcal F}_0$, so that \eqref{partialR2} is equivalent to
\begin{equation}\label{partialR3}
\lim_{n \rightarrow \infty} \frac 1 n  \left \|    S_n-{\mathbb E}(S_n| {\mathcal F}_0) - \sum_{i=1}^n D\circ T^i   \right \|_2^2 = 0 \, .
\end{equation}
Now, by Lemma \ref{D15}, \eqref{partialR3} is equivalent to 1. This concludes the proof of Proposition \ref{mainProp2}. $\square$

\subsection{The ${\mathbb L}^2$-MMA  and ${\mathbb L}^2_0$-MMA properties for H\"older functions of semi-linear processes}\label{L20MAter}

In this subsection, we first define a semi-linear process $(Y_i)_{i \in {\mathbb Z}}$ in a slightly different way. 

Let $\gamma$ in $(0, 1]$, and let ${\mathcal H}^\gamma$ be the class of a $\gamma$-H\"older functions from ${\mathbb R}$ to ${\mathbb R}$. 
Let $(\alpha_i)_{i \geq 0}$
be a sequence of functions from ${\mathcal X}$ to ${\mathbb R}$  such that
\begin{equation}\label{alphabis}
 \sum_{j \geq 0} \left \| |\alpha_j|^\gamma \right \|_{2, {\mathcal X}} < \infty\, .
\end{equation}
We then define $(A_i)_{i \geq 0}$ and  and $(Y_i)_{i \in {\mathbb Z}}$ as in \eqref{QLdef}, that is:
$A_i= \alpha_i \circ \pi_0$, 
\begin{equation}\label{QLdef2}
Y_0 =\sum_{j=0}^\infty A_j \circ T^{-j}= \sum_{j=0}^\infty \alpha_j \circ \pi_{-j} \quad \mbox{and} \quad Y_i= Y_0 \circ T^i =\sum_{j=0}^\infty A_j \circ T^{i-j} \, .
\end{equation}
Note that $Y_0$ is well defined, since 
$
\left \||Y_0|^\gamma\right \|_2 \leq \sum_{j\geq 0} \left \| |A_j |^{\gamma}\right \|_2 < \infty \, .
$
Let now 
\begin{equation}\label{defX}
X_k=f(Y_k)-{\mathbb E}(f(Y_k)) \quad \text{for $f$ in ${\mathcal H}^\gamma$}.
\end{equation}
From the previous considerations, we see that $(X_k)_{k \in {\mathbb Z}}$ is a stationary sequence of centered and square integrable random variables, adapted to the filtration 
${\mathcal F}_k= \sigma( \pi_j, j \leq k)$.

\begin{Proposition}\label{FQL}
Let $\gamma$ in $(0, 1]$, and let  $(\alpha_i)_{i \geq 0}$
be a sequence of functions from ${\mathcal X}$ to ${\mathbb R}$   satisfying \eqref{alphabis}. Then, for any sequence  $(X_k)_{k \in {\mathbb Z}}$ defined by \eqref{defX}, the property ${\mathbb L}^2$-\textup{MMA} is satisfied with $D$ defined by \eqref{Gmart}. If moreover $\sum_{j \geq 0} j \left \| |\alpha_j|^\gamma \right \|^2_{2, {\mathcal X}} < \infty$,  the property ${\mathbb L}_0^2$-\textup{MMA} is satisfied.
\end{Proposition}

\noindent{\bf Proof of Proposition \ref{FQL}.}
Let us give an upper bound for $\|P_0(X_i)\|_2$. 
Note that 
\begin{multline*}
P_0(X_i)= \\
\int  \left (f\left ( \sum_{j=0}^{i-1} \alpha_j (x_{i-j}) +A_i + \sum_{j=i+1}^\infty A_j \circ T^{i-j}\right ) - f\left ( \sum_{j=0}^{i-1} \alpha_j (x_{i-j}) +\alpha_i(x_0) + \sum_{j=i+1}^\infty A_j \circ T^{i-j}\right ) \right ) \prod_{j=0}^i \mu(d x_{i-j})\, .
\end{multline*}
Since $f \in {\mathcal H}^\gamma$, we get 
$$
|P_0(X_i)| \leq C \int |A_i-\alpha_i(x_0)|^\gamma \mu(dx_0) \, ,
$$
and consequently, since $A_i= \alpha_i \circ \pi_0$,
\begin{equation}\label{final}
\|P_0(X_i)\|_2 \leq 2 C \left \| |\alpha_i|^\gamma \right \|_{2, {\mathcal X}} \, .
\end{equation}
It follows from \eqref{final} that \eqref{H} is satisfied as soon as \eqref{alphabis} is satisfied. This proves that ${\mathbb L}^2$-MMA holds under \eqref{alphabis} with $D$ defined by \eqref{Gmart}. If moreover $\sum_{j \geq 0} j \left \| |\alpha_j|^\gamma \right \|^2_{2, {\mathcal X}} < \infty$ is satisfied, then \eqref{GL} + \eqref{H} holds, and it follows from Proposition \ref{maxProp} and Remark \ref{Rq2} that ${\mathbb L}_0^2$-MMA is satisfied. $\square$

\medskip 

\noindent {\bf Acknowledgements.} We would like to thank Christophe Cuny and Sébastien Gou\"ezel 
for useful discussions.

\end{document}